\documentclass[preprint,authoryear,12pt]{elsarticle}

\usepackage{amssymb,latexsym}
\usepackage{epsfig}
\usepackage{amsmath}
\usepackage{amssymb}
\usepackage{verbatim}
\usepackage{natbib}

\newtheorem{lem}{Lemma}
\newtheorem{thm}{Theorem}
\newtheorem{rmk}{Remark}
\journal{XXX}

\begin{document}

\begin{frontmatter}



\title{Time-varying Coefficients Estimation in Differential Equation Models with Noisy Time-varying Covariates}


\author{Heng Lian}

\address{Division of Mathematical Sciences\\School of Physical and Mathematical Sciences\\Nanyang Technological University\\Singapore, 637371\\henglian@ntu.edu.sg}

\begin{abstract}
We study the problem of estimating time-varying coefficients in ordinary differential equations. Current theory only applies to the case when the associated state variables are observed without measurement errors as presented in \cite{chenwu08b,chenwu08}. The difficulty arises from the quadratic functional of observations that one needs to deal with instead of the linear functional that appears when state variables contain no measurement errors. We derive the asymptotic bias and variance for the previously proposed two-step estimators using quadratic regression functional theory.
\end{abstract}

\begin{keyword} differential equation\sep local polynomial regression\sep measurement error\sep varying-coefficient models.



\end{keyword}

\end{frontmatter}


\section{Introduction}
Ordinary differential equations (ODEs) are widely used to describe systems in physics, chemistry, biology and medicine \citep{gardner03,cao08,miao09}. These ODEs usually involve quite a few unknown parameters that need to be estimated from observational data. Thus unlike traditional studies of dynamical systems that seek solutions for the equations, here we are concerned with the inverse problem of estimating the equations themselves given state variable measurements. Unfortunately, most ODE systems used in these applications are often complicated in form and thus do not entertain analytical solutions. Besides, the observations typically contain measurement errors and statistical methods are required to estimate these parameters. 

In general, such system can be written as
\begin{equation}\label{eqn:generalode}
\frac{d\mathbf{X}(t)}{dt}=F(\mathbf{X}(t),\boldsymbol{\beta}(t),\mathbf{V},\boldsymbol{\alpha)},
\end{equation}
where $\mathbf{X}(t)=(X_1(t),\ldots,X_p(t))^T$ are time-varying covariates, $\mathbf{V}$ are non-time-varying covariates, and $\boldsymbol{\beta}(t), \boldsymbol{\alpha}$ are time-varying and non-time-varying parameters respectively. $F$ is assumed to be known. We also assume $t\in[0,1]$ without loss of generality.
However, we do not observe $\mathbf{X}(t)$ directly. Instead we have noisy observations
\begin{equation}\label{eqn:error}
\mathbf{Y}_i=\mathbf{X}(t_i)+\boldsymbol{\epsilon}_i,
\end{equation}
where $\mathbf{Y}_i=(Y_{1i},\ldots, Y_{pi})^T$ are our actual observations and $\boldsymbol{\epsilon}_i=(\epsilon_{1i},\ldots,\epsilon_{pi})^T$ are the mean zero measurement errors assumed to be independent and identically distributed. 

Because of the importance of this problem, it has been investigated by many researchers. One approach uses classical parametric inferences such as the nonlinear least square or maximum likelihood estimation \citep{biegler86}. In this approach, optimization usually involves an iterative process, and requires using numerical methods such as Euler or Runge-Kutta. Similarly, inferences in \cite{gelman96} are based on Bayesian principle aided with Markov chain Monte Carlo methods for posterior exploration. This approach is computationally intensive since numerical
approximations to the solutions are required for each update of the parameters. 

Estimation of equation parameters that does not require numerical solutions has been proposed as early as \cite{varah82}, but seems to be largely ignored until recently. In this two-step approach, $\mathbf{X}$ and their derivatives are first estimated using a nonparametric smoother (\cite{varah82} used splines as the smoother), and in the second step the parameters in the ODEs are found based on minimizing the squared difference of the two sides of equation (\ref{eqn:generalode}) when the estimated covariates and their derivatives are plugged into the expression. This general approach is simple to implement and is taken up in some recent works \citep{chenwu08b,chenwu08,liang08,brunel08} where besides splines some of these authors used the local polynomial regression method. 

In another work, \cite{ramsay07} proposed a new method called the generalized profiling procedure. In this approach, the ODE solution is approximated by splines and both the coefficients of the basis functions and the unknown parameters in the ODEs are estimated by minimizing a penalized smoothing functional,
which reflects a trade-off between fitting the data and satisfying the ODE model. 

Both approaches described above do not required numerical solutions of ODE and have their respective advocates. Here we take the approach of the former, in particular \cite{chenwu08b,chenwu08}, and provide some new asymptotic results for a special case of (\ref{eqn:generalode}) that has not been attacked before. In particular, we consider the following ODE involving time-varying coefficients:
\begin{equation}\label{eqn:ode}
\frac{dX_1(t)}{dt}=\boldsymbol{\beta}^T(t)\boldsymbol{X}(t),
\end{equation}
where $\boldsymbol{\beta}(t)=(\beta_1(t),\ldots,\beta_p(t))^T$ are time-varying coefficients and all $X_d(t), 1\le d\le p$, are observed with measurement errors as in equation (\ref{eqn:error}).
Extension to multiple ODEs is straightforward although cumbersome in notation. We can also incorporate non-time-varying coefficients and covariates but it is regarded as simpler to analyze so we do not consider these cases. 

As far as we know, the asymptotic properties for model (\ref{eqn:ode}) are nonexistent. For the method proposed in \cite{ramsay07} and the more recent asymptotic analysis for this approach \citep{qi09}, only models involving finite-dimensional parameters are discussed. For the two-step methods, \cite{liang08,brunel08} also only consider non-time-varying parameters. \cite{chenwu08b} consider the model
\[\frac{dX(t)}{dt}=\sum_{d=1}^p \beta_d(t)Z_d(t)-g(X(t)),\]
where the functional covariates $Z_d(t)$ associated with the time-varying coefficients are observed without measurement errors and the function $g$ is known. While \cite{chenwu08} discussed a very general model
\begin{equation}\label{eqn:wugeneralode}
\frac{d\mathbf{X}(t)}{dt}=F(\mathbf{X}(t),\boldsymbol{\beta}(t))
\end{equation}
where $F$ is known, their theoretical analysis is again only limited to a very special case  
\[\frac{d\mathbf{X}(t)}{dt}=\boldsymbol{\beta}(t)-\mathbf{a}\mathbf{X}(t),\]
where the time-varying coefficients are not associated with covariates containing measurement errors and the constant $\mathbf{a}$ is known. The avoidance of these authors to analyze model (\ref{eqn:ode}) already alludes to the associated difficulties, and this is what we set out to demonstrate in this paper. 

\section{Asymptotic bias and variance}
Our problem is defined by equations (\ref{eqn:error}) and (\ref{eqn:ode}), but with the extra complication that the state variables are observed in $m$ independent experiments (say with different initial values) resulting in $m$ noisy trajectories for each state variable. More specifically, we make observations 
\[Y_{dli}=X_{dl}(t_i)+\epsilon_{dli}, 1\le d\le p, 1\le l\le m, 1\le i\le n,\]
where the state variables obey the ODEs
\[\frac{dX_{1l}(t)}{dt}=\sum_{d=1}^p\beta_d(t)X_{dl}(t), 1\le l\le m.\]
Later we will use the notations $\mathbf{Y}_{dl}=(Y_{dl1},\ldots,Y_{dln})^T$, $\boldsymbol{\epsilon}_{dl}=(\epsilon_{dl1},\ldots,\epsilon_{dln})^T$ and $\mathbf{X}_l(t)=(X_{1l}(t),\ldots, X_{pl}(t))^T$.
Note for simplicity we assume the observation times are the same for all $p$ state variables and all repeats $X_{dl}, 1\le d\le p, 1\le l\le m$. Using a two-step approach, we first estimate $X_{dl}(t)$ and the first derivative of $X_{1l}(t)$ separately using the local polynomial estimator \citep{fangijbels03}. Based on Taylor expansion, $X_{dl}(t)$ is approximated by 
\[X_{dl}(t)\approx a_0+a_1(t-t_0)+\ldots+a_q(t-t_0)^q,\]
for observation time $t$ close to a fixed point $t_0$.
Using a kernel function $K$ with a bandwidth $h$ for localization, the local polynomial estimator can be obtained by minimizing the criterion
\[\sum_{i=1}^n(Y_{dli}-\sum_{r=0}^qa_r(t_i-t_0)^r)^2K(\frac{t_i-t_0}{h}),\]
resulting in solution 
\[(T^TWT)^{-1}T^TW\mathbf{Y}_{dl},\]
where 
\[T=\left(\begin{array}{cccc}
                 1 &(t_1-t_0)&\ldots&(t_1-t_0)^q\\
		 \vdots&\vdots&&\vdots\\
	         1 &(t_n-t_0)&\ldots&(t_n-t_0)^q\\
	   \end{array}
	  \right)
\]
and $W=diag(K(\frac{t_1-t_0}{h}),\ldots,K(\frac{t_n-t_0}{h}))$. In particular, we can estimate $X_{dl}$ and its derivative $X_{dl}'=dX_{dl}/dt$ (only the derivative of $X_{1l}$ will be used though) by   
\begin{equation}\label{eqn:step11}
\hat{X}_{dl}(t_0)=\sum_{i=1}^nW_0((t_i-t_0)/h)Y_{dli},
\end{equation}
and 
\begin{equation}\label{eqn:step12}
\hat{X}'_{dl}(t_0)=\sum_{i=1}^nW_1((t_i-t_0)/h)Y_{dli},
\end{equation}
where $W_\nu(t)=e^T_{\nu,q+1}(T^TWT)^{-1}(1,ht,\ldots,h^qt^q)^TK(t),\nu=0,1$ and $e_{\nu,q+1}$ is the $(q+1)$ dimensional unit vector having $1$ as the $(\nu+1)$th component, $0$ otherwise.

In the second step, we substitute the estimates $\hat{X}_{dl}$ and $\hat{X}'_{1l}$ in the differential equation model and try to estimate the unknown coefficients $\boldsymbol{\beta}(t)=(\beta_1(t),\ldots,\beta_p(t))^T$. Again one uses local polynomial regression in this step. Around a fixed point $t_0\in (0,1)$ and approximating $\beta_d(t)$ by 
\[\beta_d(t)=\beta_{d0}+\beta_{d1}(t-t_0)+\ldots+\beta_{dq}(t-t_0)^q,\]
we obtain the local polynomial estimator $\hat{\boldsymbol{\beta}}(t)$ by minimizing the locally weighted functional
\[\sum_{i=1}^n\sum_{l=1}^m\left(\hat{X}'_{1l}(t_i)-\sum_{d=1}^p(\sum_{r=0}^q\beta_{dr}(t_i-t_0)^r)\hat{X}_{dl}(t_i)\right)^2K((t_i-t_0)/h).\]
Let 
\[Z=\left(\begin{array}{cccccc}
  \hat{X}_{11}(t_1)&\cdots&(t_1-t_0)^q\hat{X}_{11}(t_1)
			&\hat{X}_{21}(t_1)&\cdots&(t_1-t_0)^q\hat{X}_{p1}(t_1)\\
  \vdots&\vdots&\vdots&\vdots&\vdots&\vdots\\
  \hat{X}_{11}(t_n)&\cdots&(t_n-t_0)^q\hat{X}_{11}(t_n)
			&\hat{X}_{21}(t_n)&\cdots&(t_n-t_0)^q\hat{X}_{p1}(t_n)\\
 \vdots&\vdots&\vdots&\vdots&\vdots&\vdots\\
  \hat{X}_{1m}(t_1)&\cdots&(t_1-t_0)^q\hat{X}_{1m}(t_1)
			&\hat{X}_{2m}(t_1)&\cdots&(t_1-t_0)^q\hat{X}_{pm}(t_1)\\
  \vdots&\vdots&\vdots&\vdots&\vdots&\vdots\\
  \hat{X}_{1m}(t_n)&\cdots&(t_n-t_0)^q\hat{X}_{1m}(t_n)
			&\hat{X}_{2m}(t_n)&\cdots&(t_n-t_0)^q\hat{X}_{pm}(t_n)\\
      \end{array}\right)
\]
(of dimension $mn\times p(q+1)$) and let $\hat{Y}=(\hat{X}'_1(t_1),\ldots,\hat{X}'_1(t_n))^T$, the solution of the above can be written as 
\[(Z^T\mathbf{W}Z)^{-1}Z^T\mathbf{W}\hat{Y},\]
which contains estimates of $\beta_d(t_0), 1\le d\le p$ together with their derivatives, where $\mathbf{W}=I_m\otimes W=diag(W,\ldots,W)$ is the $mn\times mn$ diagonal matrix of local weights, $\otimes$ denotes the Kronecker product and $I_m$ is the $m\times m$ identity matrix. Since we are only interested in $\beta_d(t_0)$, we have the local polynomial estimator
\begin{equation}\label{eqn:sol}
\hat{\boldsymbol{\beta}}(t_0)=(I_p\otimes e^T_{0,q+1})(Z^T\mathbf{W}Z)^{-1}Z^T\mathbf{W}\hat{Y}.
\end{equation}

Note we could use different orders of polynomial and different bandwidths or even different kernels for the two steps, but we will avoid discussion on these issues since our notation is already very complicated and the results in \cite{chenwu08b} seem to suggest that these more flexible choices will not affect the asymptotic order of the estimators except for multiplicative constants for bias and variance.

We first state some standard assumptions that are used throughout the paper, which are always implicitly assumed even without mentioning. Our asymptotic results consider $m$ and $X_{dl}(\cdot)$ as fixed (or, conditional on $X_{dl}(\cdot)$) and let $n$, the number of time points, go to infinity.
\begin{enumerate}
\item[(i)] The kernel $K$ is a continuous, bounded and symmetric probability density function, with a support on $[-1,1]$.
\item[(ii)] The state variables $X_{dl}(t), 1\le d\le p, 1\le l\le m$, as well as the time-varying coefficients $\beta_d(t),1\le d\le p,$ are all three times differentiable with continuous derivatives.
\item[(iii)] The mean zero measurement errors $\epsilon_{dli}, 1\le d\le p, 1\le l\le m, 1\le i\le n$ are independent and identically distributed with finite fourth moment and its variance is denoted by $E\epsilon^2=\sigma^2$.
\item[(iv)] The observation time points $t_i, 1\le i\le n$, are independent and identically distributed with density function $f$ supported on $[0,1]$, which is continuously differentiable and bounded away from zero.
\item[(v)] The bandwidth $h$ satisfies $h\rightarrow 0$ and $nh^3\rightarrow \infty$.
\item[(vi)] Local quadratic regression is used, that is, $q=2$.
\end{enumerate}

We use several lemmas to simplify the presentation of our main results. First we have the following simple lemma concerning $Z^T\mathbf{W}Z$, which appears in (\ref{eqn:sol}). 
\begin{lem}\label{lem:zwz} $Z^T\mathbf{W}Z=nhf(t_0)[(\sum_{l=1}^m\mathbf{X}_l(t_0)\mathbf{X}_l(t_0)^T)\otimes HSH](1+o_P(1))$, where $H=diag(1,h,\ldots,h^q)$ and $S$ is a $(q+1)\times(q+1)$ matrix whose $(i,j)$ entry is $\int y^{i+j-2}K(y)dy$.
\end{lem}
\textit{Proof.} Note $Z^T\mathbf{W}Z$ can be written as
\[\sum_{i=1}^n\sum_{l=1}^m(\hat{\mathbf{X}}_l(t_i)\otimes T_i)K((t_i-t_0)/h)(\hat{\mathbf{X}}_l(t_i)\otimes T_i)^T,\]
where $T_i=(1,t_i-t_0,\ldots,(t_i-t_0)^q)^T$. Using the law of large numbers, one can show $\sum_{i=1}^n\sum_{l=1}^m(\mathbf{X}_l(t_i)\otimes T_i)K((t_i-t_0)/h)(\mathbf{X}_l(t_i)\otimes T_i)^T$ (i.e., if the covariates are observed without error) is equal to $nhf(t_0)[\sum_{l=1}^m\mathbf{X}_l(t_0)\mathbf{X}_l(t_0)^T\otimes HSH](1+o_p(1))$. The lemma easily follows from $|\hat{X}_{dl}(t)-X_{dl}(t)|=o_P(1)$.
$\Box$

The following property is well-known \citep{huangfan99,fanyao03} and is stated here only for completeness.
\begin{lem}\label{lem:wellknown}
For the weights $W_0$, $W_1$ defined immediately after (\ref{eqn:step11}) and (\ref{eqn:step12}), we have
\[\sup_{t\in[-1,1]}\sup_{t_0\in[0,1]}|nh^{\nu+1}W_\nu(t)-K_\nu(t)/f(t_0)|=o_P(1),\]
where $K_\nu(t)=e^T_{\nu,q+1}S^{-1}(1,t,\ldots,t^q)K(t), \nu=0,1$.
\end{lem}

Next we deal with the $p\times (q+1)$ dimensional vector $Z^T\mathbf{W}\hat{Y}$. First we can write 
\[Z^T\mathbf{W}\hat{Y}=\sum_{i=1}^n\sum_{l=1}^m[\hat{\mathbf{X}}_l(t_i)\otimes T_i]K(\frac{t_i-t_0}{h})\hat{X}'_{1l}(t_i).\]
A general component of this column vector is 
\[\sum_{i=1}^n\sum_{l=1}^m(t_i-t_0)^r\hat{X}_{dl}(t_i)K(\frac{t_i-t_0}{h})\hat{X}'_{1l}(t_i), 0\le r\le q, 1\le d\le p.\]
Note the appearance of $\hat{X}_{dl}(t_i)$ and $\hat{X}'_{1l}(t_i)$ together in each term of the sum is probably what deterred the researchers from studying its property.

Using (\ref{eqn:step11}) and (\ref{eqn:step12}), the above displayed expression is written as
\[\sum_{1\le i,j,k\le n}\sum_{1\le l\le m}Y_{dlj}(t_i-t_0)^rW_0(\frac{t_j-t_i}{h})W_1(\frac{t_k-t_i}{h})K(\frac{t_i-t_0}{h})Y_{1lk}=:\sum_l\mathbf{Y}_{dl}^TA_r\mathbf{Y}_{1l},\]
where the $(j,k)$ entry of the $n\times n$ matrix $A_r$, $0\le r\le q$, is defined to be
\begin{equation}\label{eqn:A}
\sum_{i=1}^n(t_i-t_0)^rW_0(\frac{t_j-t_i}{h})W_1(\frac{t_k-t_i}{h})K(\frac{t_i-t_0}{h}).
\end{equation}
The following asymptotic properties of $A_r$ are most important in deriving our main results. 
\begin{lem}\label{lem:tr}
\begin{eqnarray*}
(i) &tr(A_r)&=(C+o_P(1))h^{r-1}\\
(ii) &tr(A_r^2)&=(C+o_P(1))h^{2r-1}\\
(iii) &tr(A_rA_r^T)&=(C+o_P(1))h^{2r-1}\\
(iv) &X_{dl}^TA_rX_{1l}&=nh^{r+1}f(t_0)X_{dl}(t_0)X'_{1l}(t_0)\int y^rK(y)dy+(C+o_P(1))nh^{r+3}, 1\le d\le p\\
(v) &X_{dl}^TA_rA_r^TX_{1l}&=(C+o_P(1))h^{2r-1}+(C+o_P(1))nh^2\\
(vi) &X_{dl}^TA_r^TA_rX_{1l}&=(C+o_P(1))h^{2r-1}+(C+o_P(1))nh^2
\end{eqnarray*}
where in the above expressions, different appearances of C denote different
constants depending on the kernel $K$ and time points density $f$.
\end{lem}
\textit{Proof.} The results in the lemma are similar to those found in \cite{huangfan99}, in particular their equations (7.3), (7.6), (7.11) and (7.19). Our results are different in that we consider dense time points $t_1,\ldots, t_n$ while they consider estimation of some integral so that integrations should be replaced with summations in our case. Besides, we consider product of $W_0$ and $W_1$ in (\ref{eqn:A}) while in \cite{huangfan99} only expressions such as $W_\nu^2$ appear. Nevertheless, the calculations involved are very similar to \cite{huangfan99}, if not slightly more cumbersome. We only briefly consider the calculation of $tr(A_r)$ in the following.

Using Lemma \ref{lem:wellknown}, we can write
\begin{eqnarray*}
tr(A_r)&=&\sum_{i=1}^n\sum_{k=1}^n (t_i-t_0)^rW_0(\frac{t_k-t_i}{h})W_1(\frac{t_k-t_i}{h})K(\frac{t_i-t_0}{h})\\
&=&\sum_{i=1}^n\sum_{k=1}^n (t_i-t_0)^r \cdot \left(K_0(\frac{t_k-t_i}{h})+o_P(1)I_{\{|t_k-t_i|<h\}}\right)/(f(t_i)nh)\cdot \\
      &&\left(K_1(\frac{t_k-t_i}{h})+o_P(1)I_{\{|t_k-t_i|<h\}}\right)/(f(t_i)nh^2)\cdot K(\frac{t_i-t_0}{h})\\
&=&\sum_{i=1}^n(t_i-t_0)^r K(\frac{t_i-t_0}{h})/(f(t_i)nh^2)\int K_0(u)K_1(u)du\cdot(1+o_P(1))\\
&=&h^{r-1}\left(\int u^rK(u)du\right)\left(\int K_0(u)K_1(u)du\right)(1+o_P(1))
\end{eqnarray*}
and the result on $tr(A_r)$ is proved. One can see that the calculation strategies are quite similar to equations (7.2) and (7.3) in \cite{huangfan99}.
$\Box$

Now we can state and prove the main result in this paper.
\begin{thm} We have the following conditional bias and variance for $\hat{\boldsymbol{\beta}}(t_0)$:
\[ E(\hat{\beta}_d(t_0)-\beta_d(t_0)|t_1,\ldots,t_n)=(C_1+o_P(1))h^2+(C_2+o_P(1))\frac{1}{nh^2} \]
\[ Var(\hat{\beta}_d(t_0)|t_1,\ldots,t_n)=(C_3+o_P(1))\frac{1}{n^2h^3}+(C_4+o_P(1))\frac{1}{n}\]
for some constants $C_1, C_2, C_3$ and $ C_4$.
\end{thm}
\textit{Proof.} As observed above, a general component of $Z^T\mathbf{W}\hat{Y}$ can be written as $\sum_l\mathbf{Y}_{dl}^TA_r\mathbf{Y}_{1l}=\sum_l(\mathbf{X}^T_{dl}A_r\mathbf{X}_{1l}+\boldsymbol{\epsilon}_{dl}^TA_r\mathbf{X}_{1l}+\mathbf{X}^T_{dl}A_r\boldsymbol{\epsilon}_{1l}+\boldsymbol{\epsilon}_{dl}^TA_r\boldsymbol{\epsilon}_{1l})$, where $\mathbf{X}_{dl}=(X_{dl}(t_1),\ldots,X_{dl}(t_n))^T$ denotes the unobserved states. Using these expansions, for $d=1$, the conditional expectation of $\mathbf{Y}_{dl}^TA_r\mathbf{Y}_{1l}$ is $\mathbf{X}_{1l}^TA_r\mathbf{X}_{1l}+\sigma^2 tr(A_r)$ and the conditional variance is $4\sigma^2\mathbf{X}_{1l}^T(A_r+A_r^T)^2\mathbf{X}_{1l}+\sigma^2tr((A_r+A_r^T)^2)/2+(E\epsilon^4-3\sigma^2)\sum_{i=1}^nA_{ii,r}^2$, where $A_{ii,r}^2$ are the diagonal entries of $A_r$, while if $d\neq 1$ the conditional expectation is $\mathbf{X}_{dl}^TA_r\mathbf{X}_{1l}$ and the conditional variance is $\sigma^2[\mathbf{X}_{dl}^TA_r^TA_r\mathbf{X}_{1l}+\mathbf{X}_{dl}^TA_rA_r^T\mathbf{X}_{1l}+tr(A_rA_r^T)]$.

Based on Lemma \ref{lem:tr} and the above discussion, we can write 
\begin{equation}\label{eqn:zwy}
Z^T\mathbf{W}\hat{Y}-nhf(t_0)[\sum_{l=1}^m\mathbf{X}_l(t_0)X'_{1l}(t_0)\otimes \mathbf{w}]-\mathbf{a}_n=O_P(\mathbf{b}_n),
\end{equation}
where $\mathbf{w}=(\int K(y)dy, h\int yK(y)dy,\ldots, h^q\int y^qK(y)dy)^T$ is obtained from Lemma \ref{lem:tr} (iv), the $p\times(q+1)$ dimensional vector $\mathbf{a}_n$ is the bias term, and $\mathbf{b}_n$ is a $p\times(q+1)$ dimensional vector containing the standard deviation terms, both of which can be found from Lemma \ref{lem:tr}. The details are omitted here to avoid messy notations.

Finally, incorporating $Z^T\mathbf{W}Z$, we note

\begin{eqnarray*}
&&(I_p\otimes e^T_{0,q+1})(Z^T\mathbf{W}Z)^{-1}\left\{Z^T\mathbf{W}\hat{Y}-nhf(t_0)[(\sum_{l=1}^m\mathbf{X}_l(t_0)X'_{1l}(t_0))\otimes \mathbf{w}]-\mathbf{a}_n\right\}\\
&=&\hat{\boldsymbol{\beta}}(t_0)-(I_p\otimes e^T_{0,q+1})\frac{1}{nhf(t_0)}[(\sum_{l=1}^m\mathbf{X}_l(t_0)\mathbf{X}_l(t_0)^T)^{-1}\otimes (HSH)^{-1}]\times\\
   &&\;\;\;\left\{nhf(t_0)[\sum_{l=1}^m\mathbf{X}_l(t_0)X'_{1l}(t_0)\otimes \mathbf{w}]-\mathbf{a}_n\right\}\\
&=&\hat{\boldsymbol{\beta}}(t_0)-\boldsymbol{\beta}(t_0)-O_p((I_p\otimes e^T_{0,q+1})\mathbf{a}_n/nh)
\end{eqnarray*}

The asymptotic bias and variance is thus derived from (\ref{eqn:zwy}). $\Box$
\begin{rmk} After finding the conditional asymptotic bias and variance, it is possible, under suitable conditions, to prove asymptotic normality of $\hat{\boldsymbol{\beta}}(t_0)$, following the strategies in \cite{huangfan99}.
\end{rmk} 
\begin{rmk}
The bias and variance calculated depends on our assumptions that $X_{dl}$ is three times differentiable and local quadratic regression is used. It is possible to extend the results and get other rates when we make different assumptions on the order of smoothness of $X_{dl}$ and use local polynomial with different orders. 
\end{rmk}

\section{Conclusion}
In this paper we investigated some asymptotic properties of the two-step estimation in ODE where the time-varying coefficients are associated with noisy state variables. Asymptotic bias and variance for the estimator are found. The results presented here complement the existing results in differential equation models and make the theory more complete. The open questions include data-driven selection of the bandwidth which has not been investigated in this case and confidence interval construction. Finally, we think some extensions are possible. For example, one can use a known link function other than the identity and consider asymptotic theory for (\ref{eqn:wugeneralode}). 

\bibliographystyle{elsarticle-harv}
\bibliography{papers.txt,books.txt}







\end{document}